\documentclass{scrartcl}
\usepackage[numbers]{natbib}
\usepackage{fontenc}
\usepackage{shadethm}
\usepackage{amsthm}
\usepackage{float}
\usepackage{bm}
\usepackage{graphicx}
\usepackage{subfig}
\usepackage{stmaryrd}
\usepackage{mathrsfs} 
\usepackage{amsmath}
\usepackage{amssymb}
\usepackage{wasysym}
\usepackage{sectsty}
\usepackage[hyperfootnotes=false]{hyperref}
\usepackage[a4paper, left=2.7cm, right=2.7cm, top=2.5cm]{geometry}
\usepackage{chngcntr}
\counterwithin*{equation}{section}

\usepackage{cleveref}
\DeclareMathAlphabet{\mathpzc}{OT1}{pzc}{m}{it}

\sectionfont{\normalsize\centering}
\subsectionfont{\footnotesize\centering}

\theoremstyle{plain}
\newtheorem{thm}{Theorem}[section] 

\theoremstyle{definition}
\newtheorem{lem}[thm]{Lemma}

\newtheorem{rem}[thm]{Remark}
\newtheorem{cor}[thm]{Corollary}

\def\Xint#1{\mathchoice
	{\XXint\displaystyle\textstyle{#1}}%
	{\XXint\textstyle\scriptstyle{#1}}%
	{\XXint\scriptstyle\scriptscriptstyle{#1}}%
	{\XXint\scriptscriptstyle\scriptscriptstyle{#1}}%
	\!\int}
\def\XXint#1#2#3{{\setbox0=\hbox{$#1{#2#3}{\int}$ }
		\vcenter{\hbox{$#2#3$ }}\kern-.6\wd0}}

\def\dashint{\Xint-}
\usepackage[utf8]{inputenc}
\usepackage[german,english,russian]{babel}

\newcounter{MPequ}

\newcounter{AppA}

\begin{document}\selectlanguage{english}
\begin{center}
\normalsize \textbf{\textsf{A unique continuation theorem for exterior differential forms on Riemannian manifolds with boundary}}
\end{center}
\begin{center}
	Wadim Gerner\footnote{\textit{E-mail address:} \href{mailto:wadim.gerner@icmat.es}{wadim.gerner@icmat.es}}
\end{center}
\begin{center}
{\footnotesize	Instituto de Ciencias Matem\'{a}ticas, Consejo Superior de Investigaciones Cient\'{i}ficas, 28049 Madrid, Spain}
\end{center}
{\small \textbf{Abstract:} Aronszajn, Krzywicki and Szarski proved in \cite{AKS62} a strong unique continuation result for differential forms, satisfying a certain first order differential inequality, on Riemannian manifolds with empty boundary. The present paper extends this result to the setting of Riemannian manifold with non-empty boundary, assuming suitable boundary conditions on the differential forms. We then present some applications of this extended result. Namely, we show that the Hausdorff dimension of the zero set of harmonic Neumann and Dirichlet forms, as well as eigenfields of the curl operator (on $3$-manifolds), has codimension at least $2$. Again, these bounds were known in the setting of manifolds without boundary, so that the merit is once more the inclusion of boundary points.
\section{Introduction}
A main feature of real analytic functions, and more generally real analytic $k$-forms, defined on some connected manifold is the fact that if such a function possesses a zero of infinite order, then the whole function, respectively form, must be identically zero. This principle, namely the presence of a zero of infinite order implying being zero everywhere, is known as the strong unique continuation principle\footnote[2]{In contrast to that there is also the notion of a weak continuation principle which refers to the phenomenon that certain type of functions must be identically zero whenever they vanish on a non-empty open subset.} (sUCP). Of course, the sUCP does not hold in genreal for less regular functions, i.e. there exist smooth (but non-analytic) functions which have a zero of infinite order but do not vanish identically. But knowing that a certain type of function possess the sUCP provides additional information about its zero set. It is a standard fact that the zero set of a smooth function which has only zeros of finite order has a Hausdorff dimension of codimension at least $1$. In particular, when the function in question is not known explicitly, for instance it may arise as a solution to some elliptic PDE, it is an important step to know, in order to understand the qualitative behaviour of the solution, how the nodal sets (i.e. the zero sets) of these solutions behave. It is well-known that if the coefficients of an elliptic differential operator are analytic, then the resulting solutions are real analytic \cite[\S 3.8 \& \S 3.9]{Na85} and so they satisfy the sUCP, whereas less regular coefficients of the elliptic operator may \glqq break\grqq{} the real analyticity of the solutions. So it is of interest to know if one can preserve some features of real analyticity in this less regular setting. However, not only functions are of interest, but for instance vector fields which arise in physical applications. In these contexts one wishes to understand the flow of the vector fields, which may model particle path or electromagnetic field lines. The zero sets then make up the stationary points and are then of particular interest, since all vector fields show, at least locally, a standard behaviour around non-singular points \cite[Theorem 9.22]{L12}. So knowing the \glqq size\grqq{} of the zero set, which determines precisely the set which does \glqq not flow\grqq{} is a first step in obtaining a better understanding of the flow structure which is relevant from a physical perspective. Again, vector fields arising as solutions to certain PDE's in physics, even though they might not be real analytic, might posses the sUCP and knowing which kind of PDE leads to solutions possessing the sUCP is of relevance.
\newline
To provide some examples from physics, let us state that a classical topic in mathematical physics is to understand the behaviour of the zero sets of the Dirichlet-Laplacian\footnote[3]{Here we let the Laplace operator be a positive operator.} eigenfunctions
\[
\Delta f=\lambda f\text{, on }\Omega\subset \mathbb{R}^n\text{ with }f|_{\partial\Omega}=0,
\]
where $\Omega\subset \mathbb{R}^n$ is, let's say, a smoothly bounded domain, and to understand its complements $\Omega\setminus \{f=0\}$, see for instance \cite{HarSi89}, \cite{Log18}, \cite{SoZel11} and \cite{SoZel12} just to name a few relevant works. Of course, in the Euclidean setting, the function $f$ itself is real analytic by standard elliptic regularity results \cite[Satz 10.17]{Alt12}. However, this problem can be more generally viewed on abstract manifolds with less regular metrics. Courant's famous nodal domain theorem states for example that the number of nodal domains, i.e. the number of connected components of $\Omega\setminus\{f=0\}$ of the $k$-th Dirichlet eigenfunction is bounded above by $k$ and is at least two for $k\geq 2$, \cite[\S 5 p.365]{CouHi24}, while bounds on the $(n-1)$-dimensional Hausdorff measure of the zero set were derived in the works mentioned previously. So it is of quite some interest to understand the zero sets of Laplace eigenfunctions. Looking at less regular second order Dirichlet scalar elliptic eigenvalue problems one may ask the same questions and so it is again relevant to understand the zero set structure under less regular assumptions on the coefficients of the elliptic operator\footnote[4]{Some of the cited work indeed already deals with the less regular situation.}. Classically, the sUCP is established in the interior of the domain and this was done for instance in \cite{Mue54} for the Laplace operator, see also \cite{Aro57} for elliptic operators of second order with less regular coefficients and functions satisfying a second order differential inequality rather than necessarily an equality. These results are often based on so called Carleman estimates, see his original work \cite{Carleman39}. However, it is also of relevance to know whether or not a function has the sUCP up to the boundary, i.e. assume for the moment that a given function $f$ is smooth up to the boundary, for instance consider Dirichlet-Laplace eigenfunctions on smooth domains. Can $f$ then have a zero at the boundary, such that all the derivatives of all orders vanish at that boundary point? The new challenge in the presence of boundary is that usually when using standard elliptic estimates to establish, for example, smoothness and real analyticity of the solutions in the interior, one localises appropriately by multiplying the function by a bump function, so that if one chooses an approach via Sobolev estimates the corresponding boundary values become irrelevant. On the contrary, if we consider a boundary point, and localise the situation so that we may consider a corresponding problem on the half ball, then multiplying the solution by a bump function and integrating over the half ball, not all boundary values become irrelevant, but in fact those boundary values imposed on the disc (cutting the ball in two), i.e. the values on the boundary $\partial\Omega$, become essential. One possible approach in order to establish the sUCP at boundary points consists of the so called reflection principle, which essentially attempts to find an explicit extension of the given function from the half ball to the full ball such that this extended function still satisfies some appropriate PDE on the full ball, so that in turn one may reduce the problem to the situation of interior points. However, there is no recipe to obtain such an extension and constructing the right extension is precisely the key in extending the sUCP from interior points to boundary points, see for instance \cite{AlRoVe19} where this reflection approach was used to derive a sUCP at the boundary for a $4$-th order elliptic scalar equation. In our approach we precisely make use of this reflection principle by providing an appropriate well-behaved extension around boundary points.
\newline
\newline
The focus of our applications of the sUCP which we present in the next section will be Beltrami fields, which are eigenvector fields of the curl operator. These vector fields, also known as force-free fields, appear in the context of fluid mechanics, as solutions of the stationary Euler equations \cite[Chapter II, Remark 1.6]{AK98} and as stationary solutions of the equations of ideal magnetohydrodynamics (in the setting of constant pressure) \cite[Chapter III]{AK98}. The possible dynamics of these vector fields have been widely studied in the literature, see for instance \cite{DFGHMS86}, \cite{EP12}, \cite{EP15}, \cite{EtGr00b} and \cite{G21a} to name a few. Their zero sets were studied in detail in the real analytic setting in \cite{G21c}. Natural boundary conditions in the physical context are tangent boundary conditions, i.e. one requires the normal part of the vector fields to vanish, and a fact, known as Vainshtein's lemma, see \cite[Vainshtein's lemma]{CDGT002} for the Euclidean case, \cite[Lemma 2.1]{G21a} for the abstract manifold case and \cite{V92} for the original reference, asserts that if two curl eigenfields (of the same eigenvalue) coincide on the boundary, they must be identical everywhere. So in a sense, the behaviour of a curl eigenfield on the boundary in fact uniquely determines the behaviour of the eigenfield as a whole. In our application we strengthen this result by showing that if two eigenfields (of the same eigenvalue) coincide on the boundary on a set of Hausdorff dimension strictly larger than $1$, then they must be identical everywhere.
\newline
\newline
Let us finally point out that, as alluded in the abstract, a strong unique continuation result for Beltrami fields (more generally for abstract $k$-forms) for interior points was established in \cite{AKS62} and that the main achievement of the present paper is the extension of this result to boundary points and to apply this result in an appropriate way to Beltrami fields in order to bound the Hausdorff dimension of their boundary zero sets by $1$. Note that once a sUCP for Beltrami fields is available, it follows from \cite{B99}, in the smooth setting, that the zero set in the interior has a Hausdorff dimension of at most one. But even though we use a reflection principle to obtain a local extension of our fields which satisfies some first order elliptic equation, we loose, due to our construction, too much regularity (even if all initial quantities involved are smooth) to be able to exploit the results from \cite{B99}. Therefore, the Hausdorff dimension bound of the boundary zero set requires additional work. We derive similar Hausdorff bounds for harmonic Neumann and Dirichlet forms.
\section{Main results}
\textbf{Notation:} Throughout this paper we use the following notation: We denote by $\bar{M}$ a manifold with (possibly empty) boundary and by $M$ either its interior or, if it is specifically stated that a manifold has empty boundary, we also use $M$ to denote the full manifold to emphasise that we are dealing with a manifold with empty boundary. Further, we denote for $k\in \mathbb{N}_0$ by $\Omega^k(\bar{M})$ the maximally smooth $k$-forms on a given manifold $\bar{M}$ (maximally smooth with respect to the smoothness of the underlying manifold) and by $\mathcal{V}(\bar{M})$ the maximally smooth vector fields. Given a $C^{0,1}$-metric $g$ on $\bar{M}$ we further denote by $H^1\Omega^k(\bar{M})$ and $H^1_{\operatorname{loc}}\Omega^k(\bar{M})$ the $H^1$-regular $k$-forms and locally $H^1$-regular $k$-forms respectively. Further, given a $C^{0,1}$-metric $g$ on a $C^{1,1}$-manifold $\bar{M}$ we can decompose at a given point $p\in \partial\bar{M}$ any tangent vector field $v\in T_p\bar{M}$ as $v=v^\parallel+v^\perp$ with $v^\perp$ being the $g$-orthogonal projection of $v$ onto its normal part. The tangent part of a $k$-form $\omega\in H^1_{\operatorname{loc}}\Omega^k(\bar{M})$ is then given by $t(\omega)(p)(v_1,\dots,v_k):=\omega(p)(v^{\parallel}_1,\dots,v^{\parallel}_k)$ for $v_i\in T_p\bar{M}$, $p\in \partial\bar{M}$ and its normal part is defined as $n(\omega):=\omega-t(\omega)$ (note that the $H^1$-regularity allows us to talk in a meaningful way about traces). With this convention we set
\begin{gather}
	\nonumber
\mathcal{H}^k(\bar{M}):=\{\gamma\in H^1_{\operatorname{loc}}\Omega^k(\bar{M})|d\omega=0=\delta\omega\},\\
\nonumber
\mathcal{H}^k_N(\bar{M}):=\{\gamma\in \mathcal{H}^k(\bar{M})|n(\gamma)=0\},\\
\nonumber
\mathcal{H}^k_D(\bar{M}):=\{\gamma\in \mathcal{H}^k(\bar{M})|t(\gamma)=0\},
\end{gather}
where $d$ denotes the standard exterior derivative and $\delta$ denotes the standard co-differential. We call $\mathcal{H}^k_N(\bar{M})$ the space of harmonic Neumann $k$-forms and $\mathcal{H}^k_D(\bar{M})$ the space of harmonic Dirichlet $k$-forms. The fibre-norms induced by a given Riemannian metric $g$ are all denoted by $|\cdot|_g$. Lastly, we say that $\omega\in H^1_{\operatorname{loc}}\Omega^k(\bar{M})$ has a zero of infinite order in $1$-mean at a given point $p\in \bar{M}$, possibly a boundary point, if in one (and hence every) chart $\mu$ around $p$ each locally expressed component function $\omega_{i_1\dots i_k}$ of $\omega$ satisfies the condition
\[
\lim_{r\searrow 0}\frac{\dashint_{B_r(\mu(p))\cap \mathbb{H}^n}|\omega_{i_1\dots i_k}|dx}{r^m}=0\text{ for every fixed }m\in \mathbb{N}_0,
\]
where $\mathbb{H}^n:=\{(x^\prime,x_n)\in \mathbb{R}^{n-1}\times \mathbb{R}|x_n\geq 0\}$ denotes the upper half plane and $\dashint_A:=\frac{1}{\operatorname{vol}(A)}\int_A$ is the averaged integral. Of course, for a zero of infinite order we may omit the averaging and whenever $p$ is an interior point we may replace $\mathbb{H}^n$ by the full space $\mathbb{R}^n$.
\newline
\newline
Before we state our main result, let us recall the result from \cite{AKS62} for manifolds without boundary
\begin{thm}[AKS]
	\label{MT1}
	Let $(M,g)$ be an oriented, connected, $C^{1,1}$-smooth Riemannian $n$-manifold without boundary equipped with a $C^{0,1}$-metric $g$. Given some $k\in \mathbb{N}_0$, if $\omega\in H_{\operatorname{loc}}^1\Omega^k(M)$ satisfies
	\[
	|\delta \omega|^2_g+|d\omega|^2_g\leq C_K|\omega|^2_g\text{ a.e. on every compact subset }K\subset M,
	\]
	where $C_K>0$ is some constant which may differ for different compact sets and if $\omega$ has a zero of infinite order in $1$-mean, then $\omega=0$ a.e. on $M$.
\end{thm}
The following is our main result
\begin{thm}[Main result, strong unique continuation at boundary points]
	\label{MT2}
	Let $(\bar{M},g)$ be an oriented, connected, $C^{2,1}$-smooth Riemannian $n$-manifold with non-empty boundary equipped with a $C^{1,1}$-metric $g$. Given some $k\in \mathbb{N}_0$, if $\omega\in H_{\operatorname{loc}}^1\Omega^k(\bar{M})$ satisfies
	\begin{gather}
		\label{ME1}
	|\delta\omega|^2_g+|d\omega|^2_g\leq C_K|\omega|^2_g\text{ a.e. on every compact subset }K\subset \bar{M},
	\end{gather}
	where $C_K>0$ is a constant which may differ for different compact sets and if $\omega$ has a zero of infinite order in $1$-mean at a boundary point $p\in \partial\bar{M}$, then $\omega=0$ a.e. on $\bar{M}$, provided $\omega$ satisfies (at least) one of the following two boundary conditions
	\begin{enumerate}
		\item $t(\omega)=0$,
		\item $n(\omega)=0$.
	\end{enumerate}
\end{thm}
Let us point out first that according to \cite[Proposition 1.2.6]{S95} we have $\star t(\omega)=n(\star\omega)$ and obviously the Hodge star operator preserves the structural inequality (\ref{ME1}) as well as the property of a point of being a zero of infinite order in $1$-mean. Therefore, it is enough to prove the theorem for the boundary condition $n(\omega)=0$. Moreover, we require higher regularity in our main theorem than in the original formulation \cref{MT1}. This is because we need the normal field at the boundary to be of class $C^{1,1}$, so that in turn we may guarantee that its flow is of the same class and hence we can utilise its flow to construct an appropriate local coordinate chart of class $C^{1,1}$. In this coordinate chart the metric tensor may be locally expressed with $C^{0,1}$-coefficients and enjoys nice properties which are crucial for our construction of the reflected extension. The increased regularity is therefore necessary in our approach in order to guarantee that we may eventually apply the original result \cref{MT1} to our constructed extension. Of course, combining \cref{MT1} and \cref{MT2} one obtains a strong unique continuation property for both, interior as well as boundary, points.
\newline
\newline
As mentioned in the introduction we intend to apply this result to bound the Hausdorff dimension of the zero sets of Beltrami fields and harmonic Neumann and Dirichlet forms. We formulate these results as corollaries, even though they do not follow immediately from \cref{MT2}, but some non-trivial additional work is needed.
\begin{cor}[Hausdorff dimension of Beltrami zero sets]
	\label{MC3}
	Let $(\bar{M},g)$ be an oriented, connected, $C^{\infty}$-smooth Riemannian $3$-manifold with (possibly empty) boundary equipped with a $C^\infty$-smooth metric $g$. If $X\in \mathcal{V}(\bar{M})$ is a not identically zero, $C^\infty$-smooth vector field which is tangent to $\partial\bar{M}$ and satisfies
	\[
	\operatorname{div}(X)=0\text{ and }\operatorname{curl}(X)=\lambda X\text{ on }\bar{M}
	\]
	for a smooth function $\lambda \in C^\infty(\bar{M})$, then the Hausdorff dimension of $\{p\in \bar{M}|X(p)=0\}$ is at most $1$.
\end{cor}
Here and in the upcoming result we always compute the Hausdorff distance with respect to the natural metric on $\bar{M}$, i.e. the metric distance which is given as the smallest length of curves connecting two given points. In fact, we prove in this case not only that the Hausdorff dimension is bounded by $1$, but in fact that the zero set is $1$-$C^\infty$-smooth rectifiable in the sense that the set $\{p\in M|X(p)=0\}$ is contained in a countable union of smooth $1$-manifolds without boundary which are smoothly embedded into $M$ (into the interior of $\bar{M}$) and the set $\{p\in \partial\bar{M}|X(p)=0\}$ is contained in a countable union of smooth $1$-manifolds without boundary which are smoothly embedded in $\partial\bar{M}$. A similar, stronger, conclusion about the $(n-2)$-$C^\infty$-smooth rectifiability holds in the upcoming corollary, compare also with \cite[Corollary 3]{B99} for the statement of harmonic $k$-forms on smooth manifolds with empty boundary.
\begin{cor}[Hausdorff dimension of harmonic Neumann and Dirichlet zero sets]
	\label{MC4}
	Let $(\bar{M},g)$ be an oriented, connected, $C^{\infty}$-smooth Riemannian $n$-manifold with (possibly empty) boundary equipped with a $C^\infty$-smooth metric $g$. If $\gamma\in \mathcal{H}^k_N(\bar{M})\cup \mathcal{H}^k_D(\bar{M})$ is a not identically zero, $C^\infty$-smooth Neumann or Dirichlet form, then the Hausdorff dimension of $\{p\in \bar{M}|\gamma(p)=0\}$ is at most $(n-2)$.
\end{cor}
In both proofs of \cref{MC3} and \cref{MC4} we make use of the $C^\infty$-smoothness assumption. It would be interesting to see to what extent this regularity may be lowered.
\newline
In fact, in the real analytic setting it is standard\footnote[5]{See for example arxiv identifier: 2104.08149 Theorem A.1 for details of the proof.} that the restriction of curl eigenfields, which are tangent to the boundary, to the boundary are again real analytic vector fields, provided the boundary itself is analytic. Then it is clear, that if the Hausdorff dimension of the boundary zero set is strictly larger than one, then the restricted vector field must have a zero of infinite order and hence by real analyticity vanish on a whole boundary component. Then, in turn, the Vainshtein lemma for abstract manifolds \cite[Lemma 2.1]{G21a} implies that the original vector field must be zero everywhere. So the real strength of \cref{MC3} is that it allows us to pass from the real analytic situation to the $C^\infty$-smooth case.
\newline
Let us lastly point out that \cref{MC3} in the present form was already proven in the PhD thesis of the present author \cite[Proposition 3.2.9 (iii)]{G20Diss} and that an appropriate version of \cref{MT2} for Beltrami fields on $C^\infty$-smooth manifolds was also established \cite[Lemma 3.4.17]{G20Diss}. However, the result \cite[Lemma 3.4.17]{G20Diss} was formulated specifically for Beltrami fields, i.e. it was restricted to $3$-dimensions, to specific vector fields and a strong smoothness assumption was made. \Cref{MT2} is an extension of \cite[Lemma 3.4.17]{G20Diss} in the sense that now we allow arbitrary dimensions, we replace the restricting Beltrami field property by the much more general structural inequality (\ref{ME1}), we deal with general $k$-forms as opposed to only vector fields, i.e. $1$-forms, and we relax the regularity assumptions considerably. We also provide a new application in arbitrary dimensions in \cref{MC4}, which deals with $k$-forms.
\section{Proof of the main theorem: \Cref{MT2}}
Let us start by stating an elementary lemma
\begin{lem}
	\label{PMTL1}
	Let $\bar{M}$ be a $C^{2,1}$-smooth $n$-manifold with non-empty boundary which is equipped with a $C^{1,1}$-metric $g$. Then for every $p\in \partial\bar{M}$ there exists a $C^{1,1}$-smooth boundary coordinate chart\footnote[6]{Recall that $\mathbb{H}^n$ denotes the upper half space.} $\mu:U\rightarrow \mathbb{H}^n$ around $p$ with the following properties
	\begin{enumerate}
		\item $g_{ij}(p)=\delta_{ij}$ for all $1\leq i,j\leq n$,
		\item $\mu(p)=0$,
		\item $g_{jn}(q)=\delta_{jn}$ for all $1\leq j\leq n$ and $q\in \partial\bar{M}\cap U$.
	\end{enumerate}
\end{lem}
The proof of \cref{PMTL1} in the smooth case can be found in \cite[Lemma 3.4.19]{G20Diss} and it can be easily adapted to our situation, so that we omit the proof of \cref{PMTL1}. Let us just shortly point out that the main idea to construct these coordinates consists of considering the flowout of the inward pointing unit normal field along the boundary.
\newline
\newline
\textit{\underline{Proof of \cref{MT2}:}} Let $\omega\in H^1_{\operatorname{loc}}\Omega^k(\bar{M})$ satisfy inequality (\ref{ME1}), $n(\omega)=0$ and have a zero of infinite order in $1$-mean at some boundary point $p\in \partial\bar{M}$. By means of \cref{PMTL1} we may fix an adapted chart of class $C^{1,1}$ around $p$ and observe that the corresponding metric $g$ may be expressed in this chart with $C^{0,1}$-coefficients. By pulling back the metric and the form $\omega$ via the chart $\mu$ we obtain a $C^{0,1}$-metric and $H^1\Omega^k(B_r\cap \mathbb{H}^n)$ form, where $B_r$ is the open ball of radius $r$ centred around $0$, for some small enough $r$ which satisfies the same properties as the original form with respect to the pulled back metric and has a zero of infinite order in $1$-mean at $0$. Hence, we assume from now on that $\bar{M}=B_r\cap \mathbb{H}^n$ and that the metric $g$ is of class $C^{0,1}$ with the properties listed in \cref{PMTL1} with $p=0$. Our goal now will be to show that $\omega$ must be identically zero on $B_r\cap \mathbb{H}^n$ because then the original form has an interior zero of infinite order in $1$-mean and so the theorem will follow from \cref{MT1}. As already mentioned in the introduction we use an appropriate reflection principle. To this end we extend the metric $g$ to all of $B_r$ in the following way, where we as usual write $x=(x^\prime,x^n)\in \mathbb{R}^{n-1}\times \mathbb{R}$
\begin{gather}
	\label{PMTE1}
	\tilde{g}_{ij}(x^\prime,x^n):=\begin{cases}
		g_{ij}(x^\prime,|x^n|) & \text{ if }  i\neq n\text{ and }j\neq n\text{ or }i=j=n \\
		\operatorname{sgn}(x^n)g_{ij}(x^\prime,|x^n|) & \text{ else },
	\end{cases}
\end{gather}
where
\[
\operatorname{sgn}(z):=\begin{cases}
	+1 & \text{ if }z\geq 0 \\
	-1 & \text{ if }z<0
\end{cases}
\]
For example, we obtain for $n=3$ and $z<0$
\[
(\tilde{g}_{ij}(x,y,z))_{1\leq i,j\leq 3}=\begin{pmatrix}
	g_{11}(x,y,|z|) & g_{12}(x,y,|z|) & -g_{13}(x,y,|z|) \\
	g_{12}(x,y,|z|) & g_{22}(x,y,|z|) & -g_{23}(x,y,|z|) \\
	-g_{13}(x,y,|z|) & -g_{23}(x,y,|z|) & g_{33}(x,y,|z|)
\end{pmatrix}.
\]
We recall that by properties of our chart $\mu$ we may assume that $g_{jn}(x^\prime,0)=\delta_{jn}$ for all $(x^\prime,0)\in B_r\cap \mathbb{H}^n$ and hence the coefficients of $\tilde{g}$ may easily be seen to be Lipschitz continuous due to the Lipschitz continuity of $g$. It is obvious from definition that $\tilde{g}_{ij}$ defines a symmetric bilinear form at each point. To see that it is positive definite we employ Silvester's criterion. It is clear that the first $(n-1)$ leading principal minors are positive. Thus, we are left with observing that the determinant of the full matrix is positive at each point for small enough $r$. But this is clear because one can immediately verify that $\det{\tilde{g}}(x^\prime,x^n)=\det{g}(x^\prime,|x^n|)$ for all $(x^\prime,x^n)\in B_r$ which is positive. Hence $(\tilde{g}_{ij})_{ij}$ indeed gives rise to a well-defined $C^{0,1}$-metric on $B_r$.
\newline
As for the differential form $\omega$, we may consider its $H^1$-coefficients $\omega_{i_1\dots i_k}$ and we observe that since $\partial_n$ is a unit normal field along the boundary, the condition $n(\omega)=0$ becomes equivalent to the statement
\begin{gather}
	\label{PMTE2}
\omega_{i_1\dots i_k}(x^\prime,0)=0\text{ for all }x^\prime\text{ whenever }i_l=n\text{ for some }1\leq l\leq k.
\end{gather}
We can then define the following extension $\tilde{\omega}$ of $\omega$ to $B_r$ in terms of its coefficient functions
\[
\tilde{\omega}_{i_1\dots i_k}(x^\prime,x^n):= \begin{cases}
	\omega_{i_1\dots i_k}(x^\prime,|x^n|) &\text{ if all }i_l\neq n \\
	\operatorname{sgn}(x^n)\omega_{i_1\dots i_k}(x^\prime,|x^n|) &\text{ if at least one }i_l=n
\end{cases}
\]
With this definition it is clear that the restrictions $\tilde{w}|_{B_r\cap \mathbb{H}^n}$ and $\tilde{w}|_{B_r\cap \mathbb{H}^n_{-}}$ are of class $H^1$ on their respective domains, where we define $\mathbb{H}^n_{-}:=\{(x^\prime,x^n)\in \mathbb{R}^{n-1}\times \mathbb{R}|x^n\leq 0\}$. Further, due to the vanishing of the normal part (\ref{PMTE2}) it is also immediate that the traces of these restrictions coincide on the hyperplane $\{x^n=0\}$ and so it is standard that $\tilde{w}\in H^1(B_r)$, see also \cite[Lemma A6.9]{Alt12}.
\newline
A straightforward calculation yields
\begin{gather}
	\label{PMTE3}
	(d\tilde{\omega})_{i_1\dots i_{k+1}}(x^\prime,x^n)=\begin{cases}
		(d\omega)_{i_1\dots i_{k+1}}(x^\prime,|x^n|) & \text{ if all }i_l\neq n, \\
		\operatorname{sgn}(x^n)(d\omega)_{i_1\dots i_{k+1}}(x^\prime,|x^n|) & \text{ if some }i_l=n
	\end{cases}
\end{gather}
We then observe that $\det{\tilde{g}}(x^\prime,x^n)=\det{g}(x^\prime,|x^n|)$ on $B_r$ and
\[
\tilde{g}^{ij}(x^\prime,x^n)=\begin{cases}
	g^{ij}(x^\prime,|x^n|) &\text{ if }i\neq n\text{ and }j\neq n\text{ or }i=n=j\\
	\operatorname{sgn}(x^n) g^{ij}(x^\prime,|x^n|)&\text{ else }
\end{cases},
\]
where $\tilde{g}^{ij},g^{ij}$ as usual denote the inverse matrices. With this in mind it is also straightforward to confirm that
\begin{gather}
	\label{PMTE4}
	(d\tilde{\star}\tilde{\omega})_{i_1\dots i_{n-k+1}}(x^\prime,x^n)=\begin{cases}
		(d\star \omega)_{i_1\dots i_{n-k+1}}(x^\prime,|x^n|)& \text{ if some }i_l=n,\\
		\operatorname{sgn}(x^n)(d\star \omega)_{i_1\dots i_{n-k+1}}(x^\prime,|x^n|)& \text{ if all }i_l\neq n
	\end{cases},
\end{gather}
where $\tilde{\star}$ denotes the Hodge star operator with respect to $\tilde{g}$ and $\star$ denotes the Hodge star operator with respect to $g$. In particular, it follows from (\ref{PMTE3}) and (\ref{PMTE4}) that we control the absolute values of the coefficient functions of $d\tilde{\omega}$ and $d\tilde{\star}\tilde{\omega}$ by means of those of $d\omega$ and $d\star\omega$. Since $|\tilde{\delta}\tilde{\omega}|_{\tilde{g}}=|d\tilde{\star}\tilde{\omega}|_{\tilde{g}}$, with $\tilde{\delta}$ being the co-differential with respect to $\tilde{g}$, we conclude that $\tilde{\omega}\in H^1\Omega^k(B_r)$ satisfies inequality (\ref{ME1}) since $\omega$ does so on $B_r\cap \mathbb{H}^n$. We are left with observing that $\tilde{\omega}$ has a zero of infinite order in $1$-mean at zero simply because
\begin{gather}
	\nonumber
\int_{B_r}|\tilde{\omega}_{i_1\dots i_k}(x^\prime,x^n)|dx=\int_{B_r\cap \mathbb{H}^n}|\omega_{i_1\dots i_k}(x^\prime,x^n)|dx+\int_{B_r\cap \mathbb{H}^n_{-}}|\omega_{i_1\dots i_k}(x^\prime,-x^n)|dx
\\
\nonumber
=2\int_{B_r\cap \mathbb{H}^n}|\omega_{i_1\dots i_k}(x)|dx
\end{gather}
and so the fact that $\omega$ has a zero of infinite order in $1$-mean at zero immediately implies that so does $\tilde{\omega}$. Hence, $\tilde{\omega}$ satisfies all requirements of \cref{MT1} and therefore $\tilde{\omega}=0$ a.e. on $B_r$ so that in particular the original form $\omega$ has a zero of infinite order in $1$-mean at an interior point and a second application of \cref{MT1} to the original form $\omega$ then implies that $\omega=0$ a.e. on $\bar{M}$ (keep in mind that even though \cref{MT1} does not say anything about the boundary, the boundary is a null set in any case). $\square$
\begin{rem}
	\label{PMTR2}
	Let us point out that it is clear from (\ref{PMTE3}) that $\tilde{\omega}$ in general need not be of class $C^1$, even if the original form $\omega$ is $C^\infty$-smooth and similarly the defined metric $\tilde{g}$ in (\ref{PMTE1}) in general will be only of class $C^{0,1}$ even if the original metric $g$ is $C^\infty$-smooth.
\end{rem}
\section{Applications of \Cref{MT2}}
In this section we prove \cref{MC4} since the proof of \cref{MC3} is already contained in the present author's PhD thesis \cite[p. 390 Part III]{G20Diss}.
\newline
\newline
\textit{\underline{Proof of \cref{MC4}:}} We first note that by means of the Hodge star operator it is once more enough to consider the situation where $n(\gamma)=0$. In the $C^\infty$-smooth setting it is well-known that $\mathcal{H}^k(\bar{M})\subset C^\infty\Omega^k(\bar{M})$, i.e. all harmonic Neumann forms are smooth up to the boundary, see for instance \cite[Theorem 2.2.7]{S95}. Now we can decompose the zero set into its interior and boundary part
\[
\mathcal{Z}:=\{p\in \bar{M}|\gamma(p)=0\}=\{p\in M|\gamma(p)=0\}\sqcup \{p\in \partial\bar{M}|\gamma(p)=0\}=:\mathcal{Z}_I\sqcup\mathcal{Z}_B,
\]
where we recall that $M=\operatorname{int}(\bar{M})$ and $\sqcup$ indicates that the union is disjoint. It follows then immediately from the standard strong unique continuation result \cref{MT1} that $\gamma$ does not have any zero in $M$ of infinite order in the classical sense, i.e. if the derivatives of all orders of all coefficient functions vanish at a zero, then $\gamma$ must be identically zero. So since all zeros of $\gamma$ in $M$ are of finite order in the classical sense we conclude from \cite{B99} that the set $\mathcal{Z}_I$ is countably $(n-2)$-$C^\infty$-rectifiable. As for $\mathcal{Z}_B$ we observe that the boundary condition $n(\gamma)=0$ implies that
\[
\mathcal{Z}_B=\{p\in \partial\bar{M}|(\iota^\#\gamma)(p)=0\},
\]
where $\iota:\partial\bar{M}\rightarrow\bar{M}$ denotes the inclusion map and $\iota^\#$ its pullback. Further, let us recall \cref{PMTR2} and that the local extension which we constructed in the proof of \cref{MT2} is not of class $C^\infty$ even if all quantities involved are of class $C^\infty$. Due to this loss of regularity we cannot simply apply the result of \cite{B99} to conclude that the boundary zero set is countably $(n-2)$-rectifiable. Instead, we make use of the well-known fact that the zero set of a smooth $k$-form which has only zeros of finite order on an $(n-1)$-dimensional manifold without boundary is countably $(n-2)$-$C^\infty$-rectifiable. Hence, our goal now is to show that $\iota^\#\gamma$ has only zeros of finite order which will then conclude the proof.
\newline
\newline
We prove by induction that if $\iota^\#\gamma$ has a zero of infinite order, then $\gamma$ itself has a zero of infinite order at the boundary, so that our main theorem, \cref{MT2}, implies that $\gamma$ must be identically zero everywhere on $\bar{M}$ contradicting our assumption. We perform the induction with respect to the order of the derivatives.
\newline
\newline
\textit{Induction basis:} We note first that since $\iota^\#\gamma$ has a zero of infinite order at $p$ and since $n(\gamma)=0$, we find by working in a chart as in \cref{PMTL1} (by our regularity assumptions the corresponding chart is in fact $C^\infty$-smooth)
\[
(\partial^\beta\gamma_{i_1\dots i_k})(p)=0\text{ for all }\beta\in \mathbb{N}^{n-1}_0\times \{0\}\text{ and all }1\leq i_1,\dots,i_k\leq n.
\]
\textit{induction step:} We suppose that for some $m\in \mathbb{N}_0$ and $p\in \partial\bar{M}$
\[
(\partial^\beta \gamma_{i_1\dots i_k})(p)=0
\]
for all $\beta=(\hat{\beta},j)\in \mathbb{N}^{n-1}_0\times \mathbb{N}_0$ with $0\leq j \leq m$ and all $1\leq i_1,\dots,i_k\leq n$.
\newline
\newline
Now, let $\alpha=(\hat{\alpha},m+1)\in \mathbb{N}^{n-1}_0\times \mathbb{N}_0$. Let us consider first the situation in which $i_k=n$. We make use of the equation
\[
d\star\gamma=0,
\]
which in terms of the coefficient functions allows us to express
\[
\epsilon_{i_1\dots i_{k-1}ni_{k+1}\dots i_n}g^{i_1j_1}\dots g^{i_{k-1}j_{k-1}}g^{nj_k}(\partial_n\gamma_{j_1\dots j_{k}})
\]
in terms of the remaining coefficient functions upon which no derivative of the form $\partial_n$ is acting and where the indices $i_{k+1},\dots,i_n$ may be chosen arbitrarily. Hence, letting $\beta:=(\hat{\alpha},m)$ we may apply $\partial^\beta$ to the above expression and by means of the induction hypothesis we see that
\[
\partial^{(\hat{\alpha},0)}\left(g^{i_1j_1}\dots g^{i_{k-1}j_{k-1}}g^{nj_k}(\partial^{m+1}_n\gamma_{j_1\dots j_{k}})\right)(p)=0\text{ for all }i_1,\dots,i_{k-1}.
\]
Since $g^{nj_k}$ is constant on the boundary by choice of our coordinates we observe that we may then pull it in front of the derivative and find
\begin{gather}
	\label{AE1}
\partial^{(\hat{\alpha},0)}\left(g^{i_1j_1}\dots g^{i_{k-1}j_{k-1}}(\partial^{m+1}_n\gamma_{j_1\dots j_{k-1}n})\right)(p)=0.
\end{gather}
We now have to perform a second induction on the order of $\hat{\alpha}\in \mathbb{N}^{n-1}_0$ to see that this implies that $(\partial^\alpha\gamma_{i_1\dots i_{k-1}n})(p)=0$. Namely, we first note that the above considerations were valid for any choice of $\hat{\alpha}$ so that our induction hypothesis in particular implies, letting $\hat{\alpha}=0$ in (\ref{AE1}), $(\partial^{m+1}_n\gamma_{i_1\dots i_{k-1}n})(p)=0$. Then, if $|\hat{\alpha}|=1$, we see that all terms in (\ref{AE1}) vanish, since $(\partial^{m+1}_n\gamma_{i_1\dots i_{k-1}n})(p)=0$, except where $\partial^{(\hat{\alpha},0)}$ acts entirely on $\partial^{m+1}_n\gamma_{i_1\dots i_{k-1}n}$. But in that case we can use that $g^{i_lj_l}(p)=\delta^{i_lj_l}$ by choice of our coordinates and therefore $(\partial^{(\hat{\alpha},0)}\partial^{m+1}_n\gamma_{i_1\dots i_{k-1}n})(p)=0$ for $|\hat{\alpha}|\leq 1$ and an identical argument proves by induction that (\ref{AE1}) implies that $(\partial^{\alpha}\gamma_{i_1\dots i_{k-1}n})(p)=0$ as desired.
\newline
By multilinearity we are left with considering the situation in which all the $i_l\neq n$. To handle this case, we make use of the equation $d\gamma=0$. To this end, if we fix any set of indices $1\leq i_1<i_2<\dots<i_k<n$ it is easy to see that the condition $d\gamma=0$ implies that we can express the derivative $(\partial_n\gamma_{i_1\dots i_k})$ as a sum of the coefficient functions upon which no derivative of the form $\partial_n$ is acting. Then, as before, applying $\partial^{(\hat{\alpha},m)}$ and making use of the induction hypothesis, we conclude in this case immediately that $(\partial^\alpha\gamma_{i_1\dots i_k})(p)=0$. Thus, the induction proof is complete and we have shown that if $\iota^\#\gamma$ has a zero of infinite order at some $p\in \partial\bar{M}$, then so does $\gamma$. As explained in the beginning of the proof this shows the countably $(n-2)$-$C^\infty$-rectifiability of the boundary zero set. $\square$
\section*{Acknowledgements}
This work has received funding from the European Research Council (ERC) under the European Union’s Horizon 2020 research and innovation programme through the grant agreement 862342.
\newline
Further, this work has been partially supported by the Deutsche Forschungsgemeinschaft (DFG, German Research Foundation) – Projektnummer 320021702/GRK2326 –  Energy, Entropy, and Dissipative Dynamics (EDDy).
\bibliographystyle{plain}
\bibliography{mybibfileNOHYPERLINK}
\footnotesize
\end{document}